\documentclass[12pt]{article}
\usepackage{latexsym}
\usepackage{amssymb}
\usepackage{amsmath}
\usepackage{mathrsfs}
\usepackage[pdftex]{graphicx}
\usepackage{rawfonts}
\input{prepictex}
\input{pictex}
\input{postpictex}
\def\hok{\mbox{}\begin{picture}(10,10)\put(1,0){\line(1,0){7}}
  \put(8,0){\line(0,1){7}}\end{picture}\mbox{}}
\usepackage[OT2,OT1]{fontenc}
\def\cyr{%
\renewcommand\rmdefault{wncyr}%
\renewcommand\sfdefault{wncyss}%
\renewcommand\encodingdefault{OT2}%
\normalfont
\selectfont}
\DeclareTextFontCommand{\textcyr}{\cyr}

\def\CC{\mathbb C}

\newtheorem{main}{Theorem}

\newtheorem{thm}{Theorem}
\newtheorem{lem}{Lemma}

\newtheorem{cor}[thm]{Corollary}

\newenvironment{proof}{\medskip \noindent
{\bf Proof.}}{\hfill \rule{.5em}{1em}
\\}

\def\ZZ{{\mathbb Z}}
\def\RR{{\mathbb R}}
\def\CP{{\mathbb C \mathbb P}}
\begin{document}

\title{Einstein Metrics, 
Complex Surfaces, \\ and 
Symplectic $4$-Manifolds
}

\author{
Claude LeBrun\thanks{Supported 
in part by  NSF grant DMS-0604735.}
  }

\date{}
\maketitle

\begin{abstract}
Which  smooth compact $4$-manifolds  admit an
Einstein metric with  non-negative Einstein constant? 
A complete answer is provided in the special case of $4$-manifolds that also happen
to  admit either a complex structure or a symplectic structure. 
\end{abstract}

\bigskip

A  Riemannian manifold $(M,g)$ is said to be {\em Einstein} if it 
has constant Ricci curvature, in the sense that the function 
$$v \mapsto r (v,v)$$
on the unit tangent bundle $\{ v\in TM |~\| v\|_g=1\}$  is constant, where 
$r$ denotes the Ricci tensor of $g$. 
 This 
is of course equivalent to demanding that $g$  
  satisfy the {\em Einstein equation}  
$$r=\lambda g$$
for some real number $\lambda$. A fundamental  open problem in global  
Riemannian geometry is  to  determine  precisely
which smooth compact $n$-manifolds admit Einstein metrics. For further  background
on this problem, see  \cite{bes}. 

When  {$n=4$}, the problem  is deeply intertwined with 
 geometric and topological 
 phenomena  unique to this  dimension; 
 and our discussion here will therefore solely  focus on this 
  idiosyncratic  case. But this article   
 will focus on  even  
 narrower versions of the problem.  Let us thus first consider the special class of 
 smooth $4$-manifolds that  arise  
 from compact complex surfaces by forgetting the complex structure. 
 Which of these admit Einstein metrics? 
 If we are willing to also constrain  the 
 Einstein constant $\lambda$ to be non-negative, the following  complete answer can now be
 given:

\begin{main} \label{solo}
Let $M$ be the underlying 
 smooth $4$-manifold  
 of  a  compact complex surface.
 Then $M$ admits an  Einstein metric 
 with  $\lambda \geq 0$  if and only if  it 
is  diffeomorphic   
 to  one of  the following: a  del Pezzo surface, a $K3$ surface, 
 an Enriques surface, an Abelian surface, or a
 hyper-elliptic surface. 
\end{main}

 Recall that complex surfaces with $c_1 > 0$ are called {\em del Pezzo surfaces}. 
 The complete list of these \cite{delpezzo,cubic} consists of $\CP_1 \times \CP_1$ and 
 of $\CP_2$ blown up at $k$ points in general position, where 
 $0\leq k \leq 8$. Up to diffeomorphism, the possibilities are thus 
 $S^2 \times S^2$ and $\CP_2\# k \overline{\CP}_2$, $0\leq k \leq 8$; 
 here $\#$ denotes the connected sum, and where 
 $\overline{\CP}_2$ denotes $\CP_2$ equipped with its non-standard orientation.

 A celebrated result of   Tian \cite{tian}  asserts that 
 most del Pezzo surfaces
actually  admit $\lambda > 0$ {\em K\"ahler-Einstein} metrics; for 
earlier related  results, see   \cite{s,ty}. 
However,  there are 
  two exceptional cases  that are not covered by Tian's existence theorem; 
 namely, no K\"ahler-Einstein metric can exist on 
 on $\CP_2$ blown up at $k=1$ or $2$ points, because \cite{mats}
 both of these complex manifolds have  non-reductive
automorphism groups.  Nonetheless, an explicit  $\lambda > 0$  Einstein metric on 
$\CP_2\# \overline{\CP}_2$ was constructed  Page \cite{page}; and 
while Page's construction seemed to have nothing at all to do with 
K\"ahler geometry, Derdzi{\'n}ski \cite{derd} eventually showed
that Page's metric is in fact  {\em conformally K\"ahler} --- that is,  it
is actually a  K\"ahler metric times a smooth positive function. 
However, it  was only quite recently  \cite{chenlebweb} that 
various breakthroughs in the theory of extremal K\"ahler metrics 
made it possible to prove the existence of 
an analogous metric on $\CP_2\# 2 \overline{\CP}_2$.
The upshot is the following:

 \begin{thm}[Chen-LeBrun-Weber] \label{one}
 Let $(M^4,J)$ be any compact complex surface with $c_1 > 0$. Then 
 there is a $\lambda > 0$ Einstein metric  on $M$ 
 which  is conformally equivalent to a K\"ahler metric on $(M,J)$. 
 \end{thm}

In the  $\lambda =0$ case, the existence problem for K\"ahler-Einstein metrics was  
definitively settled  by Yau \cite{yauma}, whose  celebrated solution of the Calabi conjecture
implies   that  any compact complex manifold of 
K\"ahler type with   $c_1^\RR=0$  admits a unique Ricci-flat K\"ahler metric
in each K\"ahler class. Kodaira's classification scheme 
files the compact complex surfaces with these
 properties into four pigeon-holes \cite{bpv,GH}. First, there are the $K3$ surfaces, 
 defined as  the simply connected complex surfaces with $c_1=0$;
 they are  all deformation equivalent \cite{nicK3,kodclass}, and are thus  all diffeomorphic to 
  any smooth quartic in $\CP_3$. 
  Next,  there are the Enriques surfaces, which are $\ZZ_2$-quotients of
 suitable  $K3$ surfaces; again, there is only one diffeotype. Then  there are the Abelian surfaces, which 
 are by definition the complex tori $\CC^2/\Lambda$; obviously, these are
  all diffeomorphic to the $4$-torus. 
 And finally, there are the 
 hyper-elliptic surfaces, which are quotients of certain Abelian surfaces
 by a finite group $G$ of complex 
 affine-linear maps;  there are exactly seven possibilities for $G$, namely $\ZZ_2$, 
  $\ZZ_3$, $\ZZ_4$, $\ZZ_6$, 
 $\ZZ_2 \oplus \ZZ_2$,
 $\ZZ_2 \oplus \ZZ_4$, and 
$\ZZ_3\oplus \ZZ_3$,   
and there is exactly one diffeotype for each of these seven possibilities. 
  
The   
  existence results we have described above 
  only  produce  Einstein metrics  which are closely related to   complex structures. 
  But  what if we merely require that an Einstein metric and a  complex structure
 somehow manage 
 to uneasily  coexist  on the same manifold, without 
 necessarily being 
 on  friendly terms? Contrary to what one might expect, 
 Theorem \ref{solo} asserts
   that, provided we  constrain  the  Einstein constant to be non-negative, 
the conformally K\"ahler possibilities   already exhaust
 the entire list of possible diffeotypes.  This 
  generalizes 
  an analogous  observation regarding the $\lambda > 0$ case that was first 
  proved  in \cite{chenlebweb}.

Of course, a conformally K\"ahler,  Einstein metric
is  related not only to a  complex structure, but also to a symplectic form. 
This makes it very  tempting to look for a symplectic analog of 
Theorem \ref{solo}. What can we say, then,   about  symplectic $4$-manifolds
that  also admit  $\lambda \geq 0$ Einstein metrics?  Surprisingly enough, the answer 
turns out 
to be exactly the same!

\begin{main} \label{duet} 
Let $M$ be a smooth compact  $4$-manifold
which admits a symplectic form $\omega$. 
Then  $M$ also carries a (possibly unrelated) Einstein metric with 
$\lambda \geq 0$  if and only if  it 
is  diffeomorphic   
 to  a  del Pezzo surface, a $K3$ surface, 
 an Enriques surface, an Abelian surface, or a
 hyper-elliptic surface. 
\end{main}

What we have  learned here can thus be summarized as   follows: 

\begin{main}\label{trio} 
For a   smooth compact $4$-manifold $M$, the following statements are equivalent:
\begin{description}
\item{(i)} $M$ admits both a complex structure and an Einstein metric with 
$\lambda \geq 0$.
\item{(ii)}
 $M$ admits both a symplectic structure and an Einstein metric with 
$\lambda \geq 0$.
\item{(iii)} $M$ admits a conformally K\"ahler,  Einstein metric with 
$\lambda \geq 0$.
\end{description} 
\end{main}

One  key ingredient in  the proof of these statements is
the {\em Hitchin-Thorpe inequality} \cite{hit,tho}. Recall that the
bundle of $2$-forms over an oriented Riemannian $4$-manifold 
has an invariant decomposition 
$$\Lambda^2= \Lambda^+\oplus \Lambda^-$$
where $\Lambda^\pm$ are by definition the 
$(\pm 1)$-eigenspaces of the Hodge star operator. 
Elements of $\Lambda^+$ are called {\em self-dual $2$-forms},
and a connection on a vector bundle over $(M,g)$
is said to be {\em self-dual} if its curvature
is a bundle-valued self-dual $2$-form. This notion
is  intimately related to the $4$-dimensional 
case of the Einstein equations, because  \cite{AHS}  an oriented $4$-dimensional 
Riemannian manifold is Einstein iff the induced connection
on $\Lambda^+\to M$ is self-dual. The positivity of
instanton numbers  for solutions of the 
self-dual Yang-Mills equations 
therefore implies  the following:

\begin{lem}[Hitchin-Thorpe Inequality] \label{urht} 
Any   compact  oriented  Einstein $4$-manifold  $(M,g)$ satisfies 
 $p_1(\Lambda^+)\geq 0$, 
with  equality iff  the induced  connection on  $\Lambda^+\to M$
is flat. Moreover, the latter occurs  iff $(M,g)$ is 
finitely covered  by a Calabi-Yau $K3$ surface
or by a flat $4$-torus.  
\end{lem}
Here,  the delicate  equality case was first  cracked by  Hitchin \cite{hit}. 
An oriented Riemannian $4$-manifold $(M,g)$ induces a 
flat connection on  $\Lambda^+\to M$ iff the curvature tensor ${\mathcal R}$ of $g$ belongs to 
$\Lambda^-\otimes \Lambda^-$.  Metrics with this property
are said to be  {\em locally hyper-K\"ahler}. 
Such metrics are in particular  Ricci-flat,  so the Cheeger-Gromoll
splitting theorem \cite{bes,cg} implies  that a compact locally hyper-K\"ahler $4$-manifold
either has finite fundamental group, or else is flat. In the latter case,
Bieberbach's theorem \cite{bieber,thur3} then implies that the manifold has a finite
regular  cover which is  a flat 
$4$-torus; in the former case, it is a finite quotient of a simply connected compact manifold
with holonomy $Sp(1)$, and the choice of a parallel complex structure
then allows one to view such this universal cover as a $K3$  surface equipped
with a Ricci-flat K\"ahler metric. 

Note that, while $\Lambda^+\hookrightarrow  \Lambda^2$
depends on $g$, 
its bundle-isomophism type is  metric-independent.
In fact,  $p_1(\Lambda^+)$  actually equals the oriented homotopy invariant 
  $(2\chi + 3\tau )(M)$, where $\chi$ and $\tau$ respectively
denote the Euler characteristic and signature. For us, 
however, it is more important to notice that if $M$ admits
an orientation-compatible complex structure $J$, then 
 $$ p_1(\Lambda^+) = c_1^2 (M, J),$$
 since $\Lambda^+$ is bundle-isomorphic to $\RR\oplus K$,
 where $K=\Lambda^{2,0}$ is the canonical line bundle of $(M,J)$. 
 As a matter of convention, almost-complex structures  will henceforth always be assumed
 to be orientation-compatible. In particular, complex surfaces
 $(M,J)$ will always  be given the complex orientation, and 
 symplectic $4$-manifolds $(M,\omega)$ will always be oriented
so  that $\omega \wedge \omega$ is a volume $4$-form. 
 Thus, the Hitchin-Thorpe inequality becomes 
 \begin{equation}
 \label{enuf} 
c_1^2 (M,J) \geq 0
\end{equation}
whenever $M$ carries an almost-complex structure $J$.

To complete  the proofs of 
Theorems \ref{solo}, \ref{duet}, and \ref{trio},  let us  use 
 ($i$),  ($ii$), and  ($iii$)  to refer to the corresponding numbered statements 
 in Theorem \ref{trio}, and let us also introduce a final  numbered statement 
 \begin{description}
 \item{($iv$)} $M$ is diffeomorphic to a del Pezzo surface, a $K3$ surface, 
 an Enriques surface, an Abelian surface, or a hyper-elliptic surface. 
 \end{description}
 We have already seen that ($iv$) $\Longrightarrow$  ($iii$) $\Longrightarrow$ [($i$) and ($ii$)]. 
In light of Lemma \ref{urht} and its reformulation as  inequality   (\ref{enuf}), it  thus
suffices 
to show that
\begin{itemize}
\item if $c_1^2 > 0$, then ($i$) $\Longrightarrow$ ($ii$) $\Longrightarrow$  ($iv$); and 
\item if $c_1^2 = 0$, then ($ii$) $\Longrightarrow$ ($i$) $\Longrightarrow$  ($iv$).
\end{itemize} 

We now begin by observing that ($i$) $\Longrightarrow$ ($ii$) when $c_1^2 > 0$:

\begin{lem} Let $(M,J)$ be a compact complex surface. If $c_1^2 (M,J)>0$,
then $(M,J)$ is of K\"ahler type. In particular, $M$ admits a symplectic form $\omega$. 
\end{lem}
\begin{proof} 
Let $K=\Lambda^{2,0}$ denote the canonical line bundle  of $(M,J)$. 
Since we have assumed that 
 $c_1^2 (M)> 0$,  the Riemann-Roch theorem and Serre duality  predict that
 either $h^0(M, {\mathcal O}(K^\ell))$
or $h^0(M, {\mathcal O}(K^{-\ell}))$ must grow quadratically as $\ell \to +\infty$.
It follows \cite{bpv,kodclass} that $(M,J)$ is algebraic, and therefore  projective. 
Hence $(M,J)$ is of  K\"ahler type, as claimed. \end{proof} 

Next, we show that ($ii$) $\Longrightarrow$  ($iv$) when  $c_1^2>0$,
using  a slight 
 generalization of a result proved in   \cite{ohno}:

\begin{lem} Let $(M,\omega )$ be a symplectic $4$-manifold
with $c_1^2 > 0$. If $M$ admits a metric of non-negative scalar curvature, 
then $M$ is diffeomorphic to a del Pezzo  surface. 
\end{lem}

\begin{proof} We equip $M$ with the spin$^c$ structure
induced by any almost-complex structure adapted to the 
symplectic form $\omega$. 
Then, even  if 
$b_+(M)=1$, 
the hypothesis that  $c_1^2 >  0$ guarantees \cite{FM,spccs} 
that this spin$^c$ structure has a well-defined  Seiberg-Witten
 invariant,   
counting the  gauge-equivalence classes of 
 solutions of the {\em unperturbed}  Seiberg-Witten 
equations 
$$D_A\Phi=0, ~~F^+_A= -\frac{1}{2} \Phi \odot \bar{\Phi}$$
with multiplicities,  
for an arbitrary  Riemannian metric on $M$. 
 But since any solution of these equations must satisfy both
the Weitzenb\"ock formula 
$$0 = 2 \Delta |\Phi|^2 + 4 |\nabla \Phi|^2 + s |\Phi|^2 + |\Phi|^4$$
and the integral identity
$$\int_M|\Phi|^4 d\mu = 8\int_M|F^+_A|^2 d\mu \geq 32\pi^2 c_1^2 (M) > 0,$$
these equations have  no  solution at all  if the chosen   metric $g$ has 
scalar curvature $s\geq 0$.  Our hypotheses therefore imply that the 
Seiberg-Witten invariant must vanish for the relevant  spin$^c$ structure. However, a 
 fundamental result of Taubes \cite{taubes} asserts 
 that this invariant must be  non-zero for a symplectic $4$-manifold
 with  either  $b_+ (M)\geq 2$, or with 
$b_+(M) =1$ and $c_1\cdot [\omega ] < 0$. Our symplectic manifold therefore
 has $b_+(M) =1$ and 
$c_1\cdot [\omega ] \geq  0$.
But  since $b_+(M)=1$,  the 
 intersection form is negative-definite on the orthogonal complement
of $[\omega ]$;   our assumption that  $c_1^2 > 0$
thus implies that  $c_1\cdot [\omega ] \neq 0$,  and  our symplectic manifold therefore
has 
$b_+(M) =1$  and 
$c_1\cdot [\omega ] > 0$. 
A result of Liu \cite[Theorem B]{liu1} therefore 
tells us that $(M,\omega )$ must  be a symplectic blow-up of 
$\CP_2$ or a ruled surface. 
 Since we also have $c_1^2 > 0$, 
it follows that $M$ is diffeomorphic to $S^2 \times S^2$ or to 
$\CP_2\# k\overline{\CP}_2$ for some $k$ with $0\leq k \leq 8$. Hence 
$M$ is diffeomorphic  to a del Pezzo surface,  as claimed. 
\end{proof}

We now turn to the $c_1^2=0$ case. Recall that $b_+(M)\neq 0$ for
any symplectic $4$-manifold $(M, \omega)$, since the symplectic class
$[\omega ] \in H^2 (M, \RR)$
has positive self-intersection. 
The following observation therefore implies, in particular, 
 that ($ii$) $\Longrightarrow$ ($i$) 
when $c_1^2=0$.

\begin{lem} \label{fin} 
 Let $M$ be a smooth compact $4$-manifold 
with $p_1 (\Lambda^+) =0$ and 
 $b_+\neq 0$.  If $M$ admits an Einstein 
 metric $g$, then $g$ is Ricci-flat and  K\"ahler
 with respect to some orientation-compatible complex structure 
 $J$ on $M$. 
\end{lem}

\begin{proof} 
By Lemma \ref{urht}, the Einstein metric $g$ is locally hyper-K\"ahler.
It follows that  $g$ 
has vanishing scalar curvature $s$
and self-dual Weyl curvature $W_+$,
since these  are the trace and trace-free parts
of the $\Lambda^+\otimes \Lambda^+$ component of the 
Riemann  tensor $\mathcal R$ of $g$. 
On the other hand, since 
$b_+(M)\neq 0$, there must be a  
 self-dual harmonic $2$-form $\psi \not\equiv 0$
on $(M,g)$. 
However,  the Weitzenb\"ock formula for self-dual $2$-forms \cite{bourg,lsd}
reads
$$(d+d^*)^2 \psi = \nabla^*\nabla \psi -2 W_+(\psi, \cdot ) + \frac{s}{3}\psi,$$
so our harmonic form $\psi$ must satisfy
$$0 = \int \langle \psi ,  \nabla^*\nabla \psi \rangle d\mu =
\int |\nabla \psi |^2 d\mu$$
and we therefore have  $\nabla \psi=0$. 
In particular, the point-wise norm of $\psi$  is a non-zero constant,
and by replacing $\psi$ with a constant multiple, we may assume that 
$\|\psi\|_g \equiv\sqrt{2}$. The endomorphism 
$J : TM\to TM$ given by $v\mapsto (v\hok \psi )^\sharp$
is then  parallel, 
and satisfies $J^2=-1$. Thus $J$
 is a  complex structure   on 
$M$, and  the Ricci-flat metric $g$ now becomes 
 a K\"ahler metric on the complex surface $(M, J )$. 
\end{proof} 

Finally, we show that ($i$) $\Longrightarrow$ ($iv$) when $c_1^2=0$. 

\begin{lem} \label{beep} 
If  a smooth compact $4$-manifold $M$   admits both an
Einstein metric and a 
complex structure with $c_1^2=0$, then 
$M$ is diffeomorphic to  a $K3$ surface, an Enriques surface, 
an Abelian surface, or a hyper-elliptic surface. 
\end{lem} 
\begin{proof} By Lemma \ref{urht}, 
$M$ has a finite cover $N$ with $b_1$ even. Let $J_0$ be any  given 
complex structure on $M$, let 
$\varpi : N\to M$ denote the covering map, and let $\hat{J}_0$ denote the
pull-back of $J_0$ to $M$, so that  $\varpi$ becomes a holomorphic map from
$(N,\hat{J}_0)$ to $(M,J_0)$. 
 Since 
 the Fr\"ohlicher spectral sequence of any compact complex surface 
 degenerates at the $E_1$ level, the fact that $b_1(N)$ is even 
implies  \cite[Theorem IV.2.6]{bpv}  that 
that the real-linear  injection    $H^0(N, \Omega^1) \to H^1(N,\RR)$
 defined  by $\alpha \mapsto [ \Re e~\alpha ]$ 
is an isomorphism. Thus,  if $\varphi$ is any closed $1$-form
on $M$, $[\varpi^*\varphi ]=[ \Re e~ \alpha ]\in H^1 (N, \RR )$ 
for some $\alpha \in H^0(N, \Omega^1)$. 
 But it then  follows  that   $ [\varphi ] = [\varpi_*\varpi^* (\varphi/n) ]=[\Re e~\varpi_*(\alpha/n ) ]$,
 where $n$ is the degree of $\varpi$, and where the push-down $\varpi_*$
 is the fiber sum of the local push-forwards via the local diffeomorphism $\varpi$. 
 Hence  $H^0(M, \Omega^1) \to H^1(M,\RR)$
is also surjective, and hence an isomorphism. Thus $b_1(M)$ is even. 
Since $b_1(M) \equiv b_+(M)+1 \bmod 2$ by the integrality of the
Todd genus, it therefore follows  that $b_+(M)$ is odd, and so, in particular,  non-zero. 

Lemma \ref{fin} therefore shows that $M$ admits a complex structure $J$ 
for which the given Einstein metric $g$  is Ricci-flat and K\"ahler. 
Pulling $J$ back to the finite cover $\varpi : N\to M$
of Lemma \ref{urht} thus  realizes $(M,J)$ as the quotient of either 
a Calabi-Yau $K3$ surface or a flat Abelian surface by a 
finite group of {\em holomorphic} isometries. If the covering $\varpi$
is non-trivial, it therefore follows \cite{GH} that $(M,J)$ is either 
an Enriques surface or a hyper-elliptic surface, and  the claim therefore follows. 
 \end{proof}

By contrast, it seems much harder to determine precisely which complex
surfaces admit a general Riemannian Einstein metric if we also allow for 
the   $\lambda < 0$ case. Certainly, the Hitchin-Thorpe inequality  tells us rather immediately 
that the underlying smoth $4$-manifold of a  properly elliptic complex surface 
(that is, a surface of Kodaira dimension $1$) 
can never admit an Einstein metric.  But, by contrast,  there are  plenty 
of surfaces of general type (Kodaira dimension $2$) which {\em do} admit 
Einstein metrics. Indeed, the Aubin/Yau existence theorem \cite{aubin,yau}
tells us that there is a K\"ahler-Einstein metric with $\lambda < 0$
on any compact complex surface with $c_1 < 0$. These are precisely those 
{\em minimal} complex surfaces of general type which contain no 
$(-2)$-curves. Now, for surfaces of general type, 
 minimality turns out to have a differentiable meaning,
and not  just a holomorphic one: it means that the relevant $4$-manifold
cannot be smoothly decomposed as a connected sum $X\#\overline{\CP}_2$. 
Unfortunately, however, this is not at present known to be a 
necessary condition for the existence of an Einstein metric. 
However, one can at least prove  some weaker results
in this direction. For example \cite{lric}, if $X$ is a minimal complex surface
of general type, its $k$-point blow-up $X\# k \overline{\CP}_2$ cannot carry an 
 Einstein 
metric if $k\geq c_1^2(X)/3$. (By contrast, the Hitchin-Thorpe inequality
only gives an obstruction if $k\geq c_1^2(X)$; for an intermediate result,
see \cite{lno}.) That is, we can at least say the following:

\begin{thm} \label{five} 
Let $M$ be the underlying $4$-manifold of a compact complex surface $(M,J)$.
If $M$ admits an Einstein metric $g$, then either $M$ is as in 
Theorem \ref{solo}, or else $(M,J)$ is a surface of general type
which  is ``not too non-minimal,'' in the sense that it is  obtained from  its 
minimal model $X$ by blowing up $k <  c_1^2(X)/3$ points. 
\end{thm}

In the latter case, we of course have $c_1^2 (M) > \frac{2}{3}c_1^2(X)$.
But any minimal surface of general type satisfies \cite{bpv,GH}  the 
Noether inequality $c_1^2 (X) \geq b_+(X) -5$. Putting these together, 
and remembering that $b_+$ is unchanged by blowing up, 
we therefore obtain a non-trivial geographical inequality which, for trivial reasons,   
also happens to hold for the manifolds of Theorem \ref{solo}:

\begin{cor}
Let $M$ be the underlying $4$-manifold of a compact complex surface $(M,J)$.
If $M$ admits an Einstein metric $g$,  then $M$ is  of K\"ahler type, and 
satisfies
$$c_1^2 (M) > \frac{2}{3} \left(b_+(M) -5\right).$$
\end{cor}

In particular,  these $4$-manifolds $M$ all admit symplectic structures. 
On the other hand, there is no known result that obviously promises  such 
  a Noether-type  inequality 
for symplectic $4$-manifolds that admit   Einstein metrics. 
It would  be very interesting to prove anything    in this 
direction!

Perhaps the most fascinating  open problem in the subject is to 
 determine whether there exist Einstein metrics
on compact complex surfaces that are not 
conformally K\"ahler (with respect  to 
any  complex structure).  For surfaces with 
$c_1^2 =0$, Hitchin's  results on the boundary case of the 
Hitchin-Thorpe inequality  allow us to see  that 
no such metrics can exist. But the only other 
complex surfaces for which  such a result
has been proved  are the  
ball quotients, which saturate the Miyaoka-Yau
inequality \cite{lmo}. In a related vein, one might instead hope to improve
the ``not too non-minimal'' statement  in  Proposition 
\ref{five}. Is it really ever possible to find an Einstein metric 
on 
the underlying $4$-manifold
of a non-minimal complex surface of general type?
If so, such a metric would certainly have to 
be qualitatively different from a K\"ahler-Einstein metric, 
in many different respects!

\vfill

\noindent {\bf Dedication.} This article is dedicated to Prof.\   Akira Fujiki, and 
a preliminary version  was 
included  in the informal  lecture-note volume  
{\bf Complex Geometry in Osaka: in honour of Akira Fujiki's 60th birthday}, S. Goto {\em et al.} editors,  
Osaka University, 2008.

\vfill

\noindent 
{Department of Mathematics, SUNY, 
Stony Brook, NY 11794-3651 USA}\\
{\sc e-mail}: claude@math.sunysb.edu
 
\end{document}